\documentclass[11pt,letterpaper,draft]{article}
\usepackage[latin1]{inputenc}
\usepackage{amsmath}
\usepackage{amsfonts}
\usepackage{amssymb}
\usepackage{graphicx}
\usepackage{color}
\usepackage{comment}
\usepackage{enumerate}

\newcommand{\R}{{\ensuremath{\mathbb{R}}}}

\def\dx{{\,\mathrm{d}x}}

\def\dr{{\,\mathrm{d}r}}

\numberwithin{equation}{section}

\newtheorem{theorem}{Theorem}[section]
\newtheorem{rem}{Remark}[section]
\title{Aduen-Herron2}
\begin{document}
	\title{Nonradial solutions of  weighted elliptic superlinear problems in bounded symmetric domains}
	\author{Hugo Adu\'en\footnote{Departamento de Matem\'aticas y Estad\'\i stica, Universidad de C\'ordoba, Monter\'\i a, Colombia. E-mail address: haduen@gmail.com}, Sigifredo Herr\'on\footnote{ Escuela de Matem\'aticas, Universidad Nacional de Colombia Sede Medell\'\i n, Medell\'\i n, Colombia. E-mail address:  sherron@unal.edu.co}}
	\maketitle
	\begin{abstract}
		The present work has two objectives. First, we prove that a weight\-ed superlinear elliptic problem has infinitely many nonradial solutions in  the unit ball. Second, we obtain the same conclusion in annuli for a more general nonlinearity which also involves a weight. We use a lower  estimate of the energy level of  radial solutions with $k-1$ zeros in the interior of the domain and a simple counting. Uniqueness results due to Tanaka \cite[2008]{Tanaka1} and \cite[2007]{Tanaka2} are very useful in our approach.
	\end{abstract}
\\


\section{Introduction and statement of results}
We consider
\begin{equation}\label{theproblem}
\begin{cases}
\Delta u +K(\Vert x\Vert)\vert u\vert^{p-1}u=0, \ \text{ for } \ x\in\Omega,\\
\hfill u=0,  \text{ for }  x\in\partial\Omega, 
\end{cases}
\end{equation}
and
\begin{equation}\label{theproblem2}
\begin{cases}
\Delta u +K(\Vert x\Vert)g(u)=0, \ \text{ for } \ x\in\Omega,\\
\hfill u=0,  \text{ for }  x\in\partial\Omega, 
\end{cases}
\end{equation}
We are interested in nonradial solutions assuming $K\in C^2(\overline{\Omega})$ and positive, $\Omega$ is the unit ball for the case \eqref{theproblem} and an annulus $\Omega=\{x\in \R^N\colon a\leq \Vert x\Vert \leq b\}$  for the case \eqref{theproblem2} and $p$ is subcritical, namely $1<p<(N+2)/(N-2)$ with $N\geq 3$. 

It  is well known that some solutions of problems \eqref{theproblem} and \eqref{theproblem2} can be obtained as  critical points of the functional $J\colon H^1_0\to\mathbb{R}$ defined by
\begin{equation}\label{J}
J(u)=\int_\Omega \left( \frac{1}{2}\Vert \nabla u\Vert^2-\frac{1}{p+1}K(\Vert x\Vert)\vert u\vert ^{p+1}\right)\dx,
\end{equation}
and 
\begin{equation}\label{J_g}
J(u)=\int_\Omega \left( \frac{1}{2}\Vert \nabla u\Vert^2-K(\Vert x\Vert)G(u)\right)\dx,
\end{equation}
respectively, where $G(s)=\int_0^s g(t)\,dt.$  For simplicity, we are using the same letter $J$ in both cases. When we are looking for radial solutions to \eqref{theproblem} and \eqref{theproblem2}, the corresponding problem to be considered takes the form
\begin{equation}\label{rbvp0}
\begin{cases}
u''(r)+\frac{N-1}{r}u'(r)+K(r)\vert u(r)\vert ^{p-1}u(r)=0, \quad \text{for } r\in (0,1)\\ 
\hspace{4.34cm} u'(0)=u(1)=0,
\end{cases}
\end{equation}
and 
\begin{equation}\label{rbvp0:g}
	\begin{cases}
		u''(r)+\frac{N-1}{r}u'(r)+K(r)g(u(r))=0, \quad \text{for } r\in (a,b)\\ 
		\hspace{3.46cm} u(a)= u(b)=0,
	\end{cases}
\end{equation}
 respectively.
 From \eqref{J}, a radial solution $u$ for \eqref{theproblem} satisfies
\begin{equation}\label{J-radial}
J(u)=\left( \frac{1}{2}-\frac{1}{p+1}\right) \omega_N\int_{0}^{1}r^{N-1}K(r)\vert v(r)\vert^{p+1}\dr,
\end{equation}
where $\omega_N$ is the measure of the unit sphere in $\mathbb{R}^N$ and $v(r) = u(x)$ with $\|x\|= r$. 
In a similar fashion, if $u$ is a radial solution for \eqref{theproblem2} in the annulus $\Omega=\{x\in \R^N\colon a\leq \Vert x\Vert \leq b\}$ then,
\begin{equation}\label{J-radial:con:g}
J(u)= \omega_N\int_{a}^{b}r^{N-1}K(r)\left( \frac{g(v)v(r)}{2}-G(v)\right) \dr.
\end{equation}

From now on, all throughout the paper, $c, c_1,C, C_0,\ C_1,\ C_2,\overline{C},\ldots$ will denote generic positive constants, independent from $ u $, which may change from line to line.
\\

 In this work we prove that problems \eqref{theproblem} and \eqref{theproblem2} have infinitely many nonradial solutions in  the unit ball of $\mathbb{R}^N$ and the annuli, respectively.  For problems \eqref{theproblem} and \eqref{theproblem2}, Ramos \emph{et al} \cite{Ramos} proved the existence of a sequence $u_k$ of sign-changing  solutions whose energy levels are of order $k^\sigma$, where $\sigma=2(p+1)/(N(p-1))$, namely $J(u_k)\sim k^\sigma$. By using radial techniques, we are able to prove a lower estimate for critical levels of radial solutions $u_k$, with $k-1$ zeros and we establish that $J(u_k)\geq C(k-1)^{N\sigma}$. Then, taking into account the uniqueness results due to S. Tanaka (\cite{Tanaka1, Tanaka2}) and that the critical levels of radial solutions are more spaced, we get, by a counting argument, that most of the sign-changing solutions obtained by Ramos \emph{et al} are nonradial. 
 
 Very little about infinitely many nonradial solutions using radial tehcniques is known and, we emphasize that an upper estimate of the critical levels is not necessary. We take advantage of one result in \cite[theorem 1]{Ramos} and we complement a couple of  Tanaka's theorem by proving that \eqref{theproblem} and \eqref{theproblem2} have infinitely many nonradial solutions. Additionally, we prove that there is an infinite number of nonradial solutions considering  nonlinearities $ g(x, s) =K(\|x\|)\vert s\vert^{p-1}s$ and $ g(x, s) = K(\|x\|)g(s)$, from the   list of sign-changing solutions obtained by Ramos \emph{et al} in \cite[theorem 1]{Ramos}.
 \\

 In \cite{ACC}, an important ingredient for getting nonradial solutions was a  uniqueness result in a superlinear context.  Papers where uniqueness results have been obtained for other kinds of problem are, for example,  \cite{Tanaka3, AH}. For these, our approach does not apply. To the best of our knowledge, an estimate of critical levels as in \cite{Ramos} for sublinear problems, 
 is not known. That is why we will use the results of uniqueness due to S. Tanaka \cite{Tanaka1, Tanaka2}. Precisely, he obtained for the problem
\begin{equation}\label{rbvp}
\left\{ \begin{aligned}
u''(r) + \frac{N-1}{r}u' (r) + K(r)|u(r)|^{p-1}u(r)&=0, \quad 0<r<1,\\
u'(0)=u(1)   = 0, \ u(0) & >0,\\
&\hspace*{-6.7cm} u\ \mbox{ has exactly}\quad k-1 \  \mbox{ zeros in } (0, 1),
\end{aligned}
\right.
\end{equation}
the following result.
\begin{theorem}
Under the conditions $K\in C^2[0,1], K>0$ and
\begin{equation}\label{conditionmain}
[V(r)-p(N-2)-N+4][V(r)-p(N-2)+N]-2rV^\prime(r)<0, 
\end{equation}
where $V(r)=rK'(r)/K(r)$, the solution of problem \eqref{rbvp} exists and it is unique.
\end{theorem}
In \cite[Corollary 2.2]{Tanaka2}, S. Tanaka proved the following consequence.
\begin{theorem} Suppose $K\in C^2[a,b]$ and $K>0$. 
Assuming that:
\begin{enumerate}[\rm(H1)]
	\item $-2(N-1)\leq V(r)\leq -2$ and $V'(r)\geq 0$.
	\item The function $g$ is odd, $g\in C^1(\R)$ and $g(s)>0$ for $s>0$.
	\item $\left( g(s)/s\right)'>0$ for $s>0$.
\end{enumerate}
Then, Problem \eqref{theproblem2} has at most one radial solution $u$ with exactly $k-1$ zeros in $(a,b)$ and $u'(a)>0.$
\end{theorem}

We complement these results by proving the existence of infinitely many nonradial solutions. Our main theorems read as follows.
\begin{theorem}\label{main}
Assuming that $K\in C^2[0,1]$, $K>0$ and \eqref{conditionmain}, the problem \eqref{theproblem} has infinitely  many nonradial solutions.
\end{theorem}
\begin{theorem}\label{main2}
If $1<p<(N+2)/(N-2), \ K\in C^2[a,b], K>0$,   \textup{(H1)-(H3)} hold and, 
\begin{enumerate}
\item[\rm(H4)] There exists $C>0$ such that, for every $s>0, \  g(s) \leq C\,s^p$. 
\item [\rm(H5)]\footnote{This is the well known Ambrosetti - Rabinowitz superlinear condition.}There exists $\theta >2$ such that, for every $s>0,\ sg(s)\geq  \theta\, G(s)$.
\end{enumerate}
Then, the problem \eqref{theproblem2} has infinitely  many nonradial solutions.
\end{theorem}

\begin{rem}
As an example, the function $g(s)=\vert s\vert ^{p-1}s$ verifies \textup{(H4)} and \textup{(H5)}.	
\end{rem}

In section 2 we present some preliminaries and in section 3 we prove lower estimates of critical levels of radial solutions, which will be very important in order to prove our  theorems in section 4.

\section{Some preliminaries}
From \eqref{J-radial:con:g} and \textup{(H5)}, we observe that  
\begin{equation}\label{J-radial:con:g2}
J(u)\geq C\int_{a}^{b}r^{N-1}K(r) g(v)v(r) \dr,
\end{equation}
for every radial solution $u$ for  \eqref{theproblem2}. 
\begin{rem}\label{remark!}
Conditions \textup{(H2)-(H5)} imply:
\begin{enumerate}[\rm (a)]
\item Due to \textup{(H2)}, the function $g$ holds $sg(s)>0$ for $s\neq 0$. In addition, by using \textup{(H4)} it follows that $g(s)/s \leq C \left( sg(s)\right)^{(p-1)/(p+1)}$ for some positive constant $C$:  let us denote $\delta=(p-1)/(p+1)$. Assumption \textup{(H4)} implies
\[
	\left( \frac{g(s)}{s^p}\right) ^{1-\delta}=\left( \frac{g(s)}{s^p}\right) ^{2/(p+1)}\leq C,
\]
and thus, multiplying by $\left( \frac{g(s)}{s^p}\right) ^{\delta}$, we get
\[
	\frac{g(s)}{s^p}\leq C\left( \frac{g(s)}{s^p}\right) ^{\delta},
\]
from which the assertion follows.

\item From \textup{(H3)} and \textup{(H2)} it follows that $g'(s)>g(s)/s>0$ for $s\neq 0$. 

\item If $g(x,s):=K(\Vert x\Vert)g(s)$, assumption \textup{(H4)} implies  $g(x,s)/s\to 0$ as $s\to 0$, uniformly in $x$.

\item Again, \textup{(H4)} implies 
\[
	0\leq g(x,s)s=K(\Vert x\Vert)g(s)s\leq C\Vert K\Vert\vert s\vert^{p+1}\leq C_1(\vert s\vert^{p+1}+1).
\]

\item Because of \textup{(H5)}, 
\[
	g(x,s)s\geq \theta\, G(x,s)\geq \theta\, G(x,s)-C,
\]
where $G(x,s)=K(\Vert x\Vert)G(s)$ and $C>0$ is a constant.

\item  Hypothesis \textup{(H5)} implies that $g$ is superlinear. More exactly we have $\Big(s^{-\theta}G(s)\Big )'\geq 0$ for $s>0$ and thus, for $s>1$ we obtain $G(s) \geq G(1)s^\theta=G(1)\vert s\vert^\theta$. From this, \textup{(H2)} and \textup{(H5)} we get 
\[
	\lim\limits_{\vert s\vert \to\infty}\frac{g(s)}{s}=+\infty.
\]
\end{enumerate}
\end{rem}
In order to prove theorems \ref{main} and \ref{main2}, we shall apply theorem 1 due to  Ramos \emph{et al} in \cite{Ramos} considering special cases. In such a theorem, for problems \eqref{theproblem} and \eqref{theproblem2}, authors proved the existence of a sequence $u_k$ of sign-changing  solutions whose energy levels are of order $k^\sigma$, where $\sigma=2(p+1)/(N(p-1))$. To prove our first theorem, $\Omega$ will be the unit ball, $g(x,s)=K(\|x\|)\,|s|^{p-1}s, f(x,s)\equiv 0, \mu=p+1$  and we choose any number 
\[
\nu\in\left( 0,\frac{N+2-p(N-2)}{2}\right), 
\]
in order to obtain condition (1.4) in \cite{Ramos}. In this context, such a theorem is established as follows.
\begin{theorem}\label{ramos}
Assuming that $N\geq3, \ 1 <p< (N+2)/(N-2)$,

the problem 
\[
\Delta u +K(\Vert x\Vert)\vert u\vert^{p-1}u=0; \, u\in H_0^1(\Omega),
\]
 admits a sequence of sign-changing solutions $(u_k)_{k\in\mathbb{N}}$ whose energy levels $J(u_k)$ satisfy
\begin{equation}\label{ineq1ramos}
c_1k^{\sigma}\leq J(u_k)\leq c_2k^{\sigma},
\end{equation}
for some $c_1, c_2>0$ with $\sigma=\frac{2(p+1)}{N(p-1)}$.
\end{theorem}
To prove our second theorem, we will take $\Omega$ as an annulus, $g(x,s)=K(\|x\|)\,g(s)$, $f(x,s)\equiv 0, \mu=\theta$ and we choose
\[
\nu\in\left( 0,\frac{\theta(N+2-p(N-2))}{2(p+1)}\right), 
\]
with the aim that condition (1.4) in \cite{Ramos} holds; further, the above  remarks imply all conditions in \cite[Theorem 1]{Ramos} are satisfied and hence, its conclusion give us a sequence of sign-changing solutions $(u_k)_{k\in\mathbb{N}}$ whose energy levels $J(u_k)$ satisfy \eqref{ineq1ramos}. 	

\section{Lower estimates of critical leves}

In this section we obtain estimates of the critical levels corresponding to a radial solution $u_k$ with $k-1$ zeros for the problems \eqref{theproblem} and \eqref{theproblem2}. More exactly, in order to prove our first main result we establish an estimate from below of $J(u_k)$ where $u_k$ is a radial solution of \eqref{rbvp} with $k-1$ zeros in $(0,1)$. Then, the same estimate will be gotten for a radial solution $u_k$ with $k-1$ zeros for the problem \eqref{theproblem2}.
\begin{theorem}\label{boundbelow}
	Let $\delta\colon=(p-1)/(p+1)$. There exists a constant $C>0$  such that for all solution $u_k\equiv u$ of \eqref{rbvp}, we have 
	\begin{equation}\label{ineq2Kaji}
	J(u)\ge C(k-1)^{N\sigma}.
	\end{equation}
\end{theorem}

\begin{theorem}\label{boundbelow}
	 There exists a constant $C>0$  such that for all solution $u_k\equiv u$ with $k-1$ zeros of  \eqref{rbvp0:g}, we have 
	\begin{equation}\label{ineq2Kaji}
	J(u)\ge C(k-1)^{N\sigma}.
	\end{equation}
\end{theorem}

\section{Proof of theorems \ref{main} and \ref{main2}}

By using theorems of section 3 and a counting argument, we can show our main results.

\vspace{.5cm}

\bibliographystyle{amsplain}

\begin{thebibliography}{10}
	\bibitem{AC} Adu\'en, H. and Castro, A. \emph{Infinitely Many Nonradial Solutions to a Superlinear Dirichlet Problem}.  Proceedings of the American Mathematical Society, vol. 131, No. 3, 2003, pp. 835-843.
	\bibitem{ACC} Adu\'en, H.; Castro, A.; Cossio, J.\emph{
	Uniqueness of large radial solutions and existence of nonradial solutions for a superlinear Dirichlet problem in annuli}. J. Math. Anal. Appl., Volume 337, Issue 1, 2008, pp. 348-359.


\bibitem{AH}  Adu\'en, H. and  Herr\'on, S. \emph{On the uniqueness of sign-changing solutions to a semipositone
problem in annuli}, Rev. Integr. Temas Mat. 34 (2016), No. 2, 207-224.

\bibitem{Amadori}  Amadori, A. L. and  Gladiali,
F. \emph{Nonradial sign changing solutions to Lane-Emden problem in an annulus}. Nonlinear Analysis,
Volume 155, 2017, pp. 294-305.
	
	\bibitem{Bartsch} Bartsch, T. and Willem,
M. \emph{Infinitely Many Nonradial Solutions of a Euclidean Scalar Field Equation}. Journal of Functional Analysis,
Volume 117, Issue 2,
 	1993, 	Pages 447-460.
 	
	
	\bibitem {Hartman}  Hartman, P. \emph{Ordinary Differential Equations: Second Edition}. Volumen 38  Classics in Applied Mathematics, SIAM. 2002.
	\bibitem{hirano} Hirano, N. and Mizoguchi, N. \emph{Nonradial solutions of semilinear elliptic equations on annuli}. J. Math. Soc. Japan, No 1, (1994), pp. 111-117.
	
	
	\bibitem{Kajikiya} Kajikiya, R. \emph{Sobolev norms of radially symmetric oscillatory solutions for superlinear elliptic equations}. Hiroshima Math. J. 20 (1990), No. 2, 259-276.
	\bibitem{Kajikiya2} Kajikiya, R. \emph{Non-radial solutions with orthogonal subgroup invariance for semilinear Dirichlet problems}. Topol. Methods Nonlinear Anal. 21 (2003), No. 1, 41-51.
		
	\bibitem {Ramos}   Ramos, M.; Tavares, H.; Zou, W. \emph{A Bahri-Lions theorem revisited.} Adv. Math. 222, No. 6, 2173-2195 (2009).
	
	\bibitem{Rolando} Rolando, S.  \emph{Multiple nonradial solutions for a nonlinear elliptic problem with singular and decaying radial potential}. Adv. Nonlinear Anal. 8:  (2017), pp. 885-901.
	
	\bibitem{Serra}Serra, E. \emph{Non-radial positive solutions for the Hénon equation with critical growth}.
	Cal. Var. Partial Differential Equations 23 (2005), No. 3, 301-326.
	
\bibitem {Tanaka1}   Tanaka, S. \emph{Uniqueness of nodal radial solutions of superlinear elliptic equations in a ball.} Proc. Roy. Soc. Edinburgh Sect. A 138 (2008), No. 6, p. 1331-1343.

\bibitem {Tanaka2}   Tanaka, S. \emph{On the uniqueness of solutions with prescribed numbers of zeros for a two-point boundary value problem.} Differential Integral Equations Volume 20, Number 1 (2007), 93-104.

\bibitem{Tanaka3} Tanaka,
S. \emph{Uniqueness and nonuniqueness of nodal radial solutions of sublinear elliptic equations in a ball},
Nonlinear Analysis: Theory, Methods \& Applications,
Volume 71, Issue 11,
2009,
5256-5267.

\bibitem{Wei} Wei, J. and  Yan, S. \emph{Infinitely many nonradial solutions for the
	H\'enon equation with critical growth}. Rev. Mat. Iberoam. 29 (2013), No. 3, 997-1020.

\bibitem{Zhang} Zhang, H. and  Zhang, F. \emph{Infinitely many radial and nonradial solutions for a Choquard equation with general nonlinearity}. Applied Mathematics Letters,
 Volume 102,
 2020, 106142.

\end{thebibliography}

\end{document}